*Classical Injective Solutions in the Large in Incompressible Nonlinear Elasticity*

T.J. HEALEY
Department of Mathematics
Cornell University
e-mail: tjh10@cornell.edu

**Abstract**

We consider a general class of parametrized displacement boundary value problems in incompressible nonlinear elasticity. We prove the existence of an unbounded solution branch of classical injective solutions emanating from the unforced stress-free reference configuration – a sharp global implicit function theorem. The ramifications of this are: There exists at least one solution for any given applied loading, and/or there exists solutions of arbitrarily large norm for some finite loading. We refer to this as *solutions in the large.* The nonlinear constraint equation enforcing incompressibility obviates the direct use of either the Leray-Schauder degree or an oriented degree previously developed for compressible problems. Instead we employ the nonlinear Fredholm degree for proper maps constructed in [5], [21]. We establish the requisite admissibility properties in the presence of the constraint equation before proving the main results.

**1. Introduction**

The goal of this work is the same as that of [9] and [11], but in the context of incompressible nonlinear elasticity. Global branches of classical, injective solutions were obtained in those works for a general class of parametrized displacement boundary value problems in compressible nonlinearly elasticity. The main result of [9] is the existence of an unbounded solution branch emanating from the unforced, stress-free reference configuration - a sharp global implicit function theorem. We emphasize the ramifications: There exist solutions of arbitrarily large norm for some finite loading, and/or there exists at least one solution for any given applied loading, the latter as in the classical result of Leray and Schauder [18]. We refer to this alternative as *solutions in the large*, which we rigorously establish here in the incompressible case. The starting point in [9] is the application of an oriented, $C^2$ nonlinear Fredholm degree developed in [10], while in [11] the Leray-Schauder degree is shown to be applicable to the same class of displacement problems considered in [9].

The nonlinear constraint equation enforcing incompressibility (or local volume preservation) obviates the direct use of either of the above-mentioned degrees in the class of problems considered here. On the other hand, incompressibility eliminates a real source of difficulties in the analysis of compressible problems, viz., the unbounded growth of the stored-energy density as the local volume ratio approaches zero. As such, we arrive at an elliptic system posed on an entire Banach space. This setup is tailor-made for the nonlinear Fredholm degree of Fitzpatrick, Pejsachowicz and Rabier [5], [21], which we employ to obtain our results. That degree (henceforth referred to as the FPR degree) is based on the parity of curves in simply-connected open subsets of a Banach space. We note that local solutions for this class of incompressible problems have been obtained in [17] via the implicit function theorem. Our results here show that such paths are only the "start" of global solution continua.

An outline of the paper is as follows. In Section 2 we present our formulation. Assuming strong ellipticity compatible with incompressibility (cf. [6], [27]), we consider a general class of parametrized problems characterized by displacement-dependent or "live" body-force densities and placement conditions on the boundary. The incompressibility constraint is treated as a field equation, while the pressure appears naturally in the equilibrium equations. We choose to work in spaces of Hölder-continuous functions, which is particularly convenient for the treatment live body forces. In Section 3 we obtain the existence of local solution paths via the implicit function theorem. In particular, we give the details for establishing the Fredholm property (index zero) via the elliptic estimates of Agmon, Douglis and Nirenberg [1].

In Section 4 we state a general global result for solution continua with a Rabinowitz-type alternative, cf. [23]. We then prove two key properties of our mapping that are needed for the application of the abstract FPR degree: (i) the nonlinear Fredholm property; (ii) properness. The first step (i) is similar to that given in Section 3, based on the elliptic estimates and the stability of the Fredholm index. Due to the presence of the constraint equation, the governing elliptic system is not quasi-linear. Hence, our proof of properness (ii) is not routine. Next we summarize the FPR degree and then use it to prove the theorem.

We make two additional, physically reasonable hypotheses in Section 5: The stored energy is strictly quasiconvex at the stress-free reference configuration (only), and the domain is star-shaped. A result of [16] then implies that the unloaded reference configuration corresponds to the unique solution of the problem in the absence of loading. Hence, we deduce that the solution continuum established in Section 4 is unbounded, i.e., there exist solutions in the large. We give some final remarks in Section 6.

*Notation*

For Banach spaces $X$ and $Y$, $L(X,Y)$ denotes the space of all bounded linear transformations of $X$ into $Y$, while $L(X) := L(X,X)$. Vectors in $\mathbb{R}^3$ and tensors in $L(\mathbb{R}^3)$ are denoted by bold lower-case and upper-case symbols, respectively. The Euclidean inner-product of two vectors in $\mathbb{R}^3$ is denoted $\mathbf{u} \cdot \mathbf{v}$, with $|\mathbf{u}| = \sqrt{\mathbf{u} \cdot \mathbf{u}}$. We generally write $A[x] \in Y$ for the value of $A \in L(X,Y)$ at $x \in X$. For the special case $\mathbf{A} \in L(\mathbb{R}^3)$, we simply write $\mathbf{Ax}$ instead of $\mathbf{A}[\mathbf{x}]$; $\mathbf{I} \in L(\mathbb{R}^3)$ denotes the identity. We define

$$GL^+(\mathbb{R}^3) = \{\mathbf{A} \in L(\mathbb{R}^3) : \det \mathbf{A} > 0\},$$
$$U(\mathbb{R}^3) = \{\mathbf{A} \in L(\mathbb{R}^3) : \det \mathbf{A} = 1\},$$
$$SO(3) = \{\mathbf{A} \in GL^+(\mathbb{R}^3) : \mathbf{A}^{-1} = \mathbf{A}^T\},$$
$$TS(\mathbb{R}^3) = \{\mathbf{A} \in L(\mathbb{R}^3) : \mathbf{A}^T = \mathbf{A}, \ tr\mathbf{A} = \mathbf{0}\}.$$

The cofactor tensor of $\mathbf{F}$ is denoted $Cof \ \mathbf{F}$; $Cof \ \mathbf{F} = (\det \mathbf{F})\mathbf{F}^{-T}$ for $\mathbf{F} \in GL^+(\mathbb{R}^3)$. Also, $\mathbf{A} \cdot \mathbf{B} = tr(\mathbf{AB}^T)$ denotes the inner-product on $L(\mathbb{R}^3)$, and $|\mathbf{A}| = \sqrt{\mathbf{A} \cdot \mathbf{A}}$. Finally, $\mathcal{C}[\mathbf{A}] \in L(\mathbb{R}^3)$ is the value of the fourth-order tensor $\mathcal{C} \in L(L(\mathbb{R}^3))$ at $\mathbf{A} \in L(\mathbb{R}^3)$.

## 2. Formulation

Throughout this work we assume $\Omega \subset \mathbb{R}^3$ is a bounded domain, with a $C^3$ boundary $\partial \Omega$, which we take as our reference configuration. The outward unit normal at $\mathbf{x} \in \partial \Omega$ is denoted $\mathbf{n}(\mathbf{x})$. We denote a deformation by $\mathbf{f} : \bar{\Omega} \to \mathbb{R}^3$, i.e., $\mathbf{y} = \mathbf{f}(\mathbf{x})$ is the position of the material point in the deformed configuration, $\mathbf{f}(\bar{\Omega})$, that occupies $\mathbf{x} \in \bar{\Omega}$. For a smooth vector field, $\mathbf{v} : \bar{\Omega} \to \mathbb{R}^3$, the gradient or total derivative is the second-order tensor-valued function $\nabla \mathbf{v} : \bar{\Omega} \to L(\mathbb{R}^3)$, and the divergence is given by $\nabla \cdot \mathbf{v} = tr(\nabla \mathbf{v})$. The divergence of a smooth tensor field, $\mathbf{S} : \bar{\Omega} \to L(\mathbb{R}^3)$, is the vector field $\nabla \cdot \mathbf{S}$ defined by $(\nabla \cdot \mathbf{S}) \cdot \mathbf{a} = \nabla \cdot (\mathbf{S}^T \mathbf{a}) \; \forall \; \mathbf{a} \in \mathbb{R}^3$. We let $\mathbf{F} := \nabla \mathbf{f}$ denote the deformation gradient; incompressibility is expressed via the constraint $\mathbf{F} \in U(\mathbb{R}^3)$ in $\bar{\Omega}$.

We assume that $\bar{\Omega}$ is occupied by an incompressible, homogeneous, hyperelastic body, i.e., there is a sufficiently smooth stored-energy density $W : U(\mathbb{R}^3) \to \mathbb{R}$ such that $W_{\mathbf{F}} := dW / d\mathbf{F}$ represents to the constitutively determined part of the nominal or Piola-Kirchhoff stress tensor. The stored energy function $W(\cdot)$ possesses a smooth extension to $GL^+(\mathbb{R}^3)$ (cf. [6], [16]), which we henceforth presume, that conforms with material objectivity:

$$W(\mathbf{QF}) = W(\mathbf{F}) \; \forall \; \mathbf{Q} \in SO(3). \tag{2.1}$$

We assume that the reference configuration is stress free, viz.,

$$W_{\mathbf{F}}(\mathbf{I}) = 0, \tag{2.2}$$

and we suppose that the body is subjected to one-parameter families of body-force densities, $\tilde{\mathbf{b}}(\lambda, \nabla \mathbf{f}, \mathbf{f}, \cdot)$ in $\Omega$, and homogeneous deformation fields $\mathbf{f}(\mathbf{x}) = \mathbf{A}(\lambda)\mathbf{x}$ on $\partial \Omega$, $\mathbf{A}(\lambda) \in U$, where $\lambda \in \mathbb{R}$ is a parameter. Let $\mathbf{S}(\mathbf{x})$ denote the total nominal stress tensor at $\mathbf{x}$. For $\mathbf{F} \in U(\mathbb{R}^3)$, we have

$$\mathbf{S} = W_{\mathbf{F}}(\mathbf{F}) - p Cof \, \mathbf{F}, \tag{2.3}$$

where $p$ denotes the constitutively indeterminate pressure (enforcing incompressibility). Recall that $\nabla \cdot (Cof \, \nabla \mathbf{v}) \equiv \mathbf{0}$ for any smooth vector field on $\bar{\Omega}$. In view of (2.3), the well-known local form of the equilibrium equations then reads

$$\nabla \cdot \mathbf{S} = \nabla \cdot (W_{\mathbf{F}}(\nabla \mathbf{f})) - [Cof \, \nabla \mathbf{f}] \nabla p = -\tilde{\mathbf{b}}. \tag{2.4}$$

The fourth-order elasticity tensor at $\mathbf{F} \in U(\mathbb{R}^3)$, defined by

$$\mathcal{C}(\mathbf{F}) := \frac{d^2 W(\mathbf{F})}{d\mathbf{F}^2}, \tag{2.5}$$

is assumed to satisfy the strong ellipticity condition

$$(\mathbf{a} \otimes \mathbf{c}) \cdot \mathcal{C}(\mathbf{F})[\mathbf{a} \otimes \mathbf{c}] \geq \gamma |\mathbf{a}|^2 |\mathbf{c}|^2 \; \forall \mathbf{F} \in U(\mathbb{R}^3), \mathbf{a}, \mathbf{c} \in \mathbb{R}^3 \text{ s.t. } \mathbf{a} \cdot [Cof \mathbf{F}]\mathbf{c} = 0, \tag{2.6}$$

where $\gamma > 0$ is a constant.

Incorporating (2.5) into (2.4), and accounting for the incompressibility constraint, we arrive at the field equations:

$$\begin{aligned} \mathcal{C}(\nabla \mathbf{f}) \nabla^2 \mathbf{f} - [Cof(\nabla \mathbf{f})] \nabla p + \tilde{\mathbf{b}}(\lambda, \nabla \mathbf{f}, \mathbf{f}) &= \mathbf{0} \text{ in } \Omega, \\ \det(\nabla \mathbf{f}) &= 1 \text{ in } \Omega, \\ \mathbf{f}(\mathbf{x}) &= \mathbf{A}(\lambda)\mathbf{x} \; \mathbf{x} \in \partial \Omega. \end{aligned} \tag{2.7}$$

In components relative to the standard basis, the first term on left side of (2.7)$_1$ reads

$$(\mathcal{C}(\nabla \mathbf{f})\nabla^2 \mathbf{f})_i = \mathcal{C}_{ijkl}(\nabla \mathbf{f})\frac{\partial^2 f_k}{\partial x_j \partial x_l}. \tag{2.8}$$

We make the substitution $\mathbf{f}(\mathbf{x}) = \mathbf{A}(\lambda)\mathbf{x} + \mathbf{u}(\mathbf{x})$ in (2.7), to obtain the final form of the field equations:

$$\begin{aligned}
\mathcal{C}(\mathbf{A}(\lambda) + \nabla \mathbf{u})\nabla^2 \mathbf{u} - [Cof(\mathbf{A}(\lambda) + \nabla \mathbf{u})]\nabla p + \mathbf{b}(\lambda, \nabla \mathbf{u}, \mathbf{u}) &= \mathbf{0} \text{ in } \Omega, \\
\det(\mathbf{A}(\lambda) + \nabla \mathbf{u}) - 1 &= 0 \text{ in } \Omega, \\
\mathbf{u} &= \mathbf{0} \text{ on } \partial\Omega,
\end{aligned} \tag{2.9}$$

where $\mathbf{b}(\lambda, \nabla \mathbf{u}, \mathbf{u}, \mathbf{x}) := \tilde{\mathbf{b}}(\lambda, \mathbf{A}(\lambda) + \nabla \mathbf{u}, \mathbf{A}(\lambda)\mathbf{x} + \mathbf{u}, \mathbf{x})$.

We assume the following smoothness conditions:

$$\begin{aligned}
W &\in C^4(GL^+(\mathbb{R}^3)), \\
\mathbf{b} &\in C^3(\mathbb{R} \times L(\mathbb{R}^3) \times \mathbb{R}^3 \times \overline{\Omega}, \mathbb{R}^3), \\
\mathbf{A}(\cdot) &\in C^1(\mathbb{R}, U(\mathbb{R}^3)).
\end{aligned} \tag{2.10}$$

In addition, we assume zero loading at $\lambda = 0$:

$$\begin{aligned}
\mathbf{b}(0, \cdot) &\equiv \mathbf{0} \text{ on } \overline{\Omega}, \\
\mathbf{A}(0) &= \mathbf{I}.
\end{aligned} \tag{2.11}$$

Next we express (2.7) in abstract operator form. Let $X$ and $Y$ denote the Banach spaces defined below:

$$\begin{aligned}
X &:= \{(\mathbf{u}, p) \in C^{2,\alpha}(\overline{\Omega}, \mathbb{R}^3) \times C_a^{1,\alpha}(\overline{\Omega}) : \mathbf{u} = \mathbf{0} \text{ on } \partial\Omega\}; \quad \|(\mathbf{u}, p)\|_X = \|\mathbf{u}\|_{2,\alpha;\overline{\Omega}} + \|p\|_{1,\alpha;\overline{\Omega}}, \\
Y &:= \{(\boldsymbol{\tau}, \rho) \in C^{\alpha}(\overline{\Omega}, \mathbb{R}^3) \times C_a^{1,\alpha}(\overline{\Omega})\}; \quad \|(\boldsymbol{\tau}, \rho)\|_Y = \|\boldsymbol{\tau}\|_{\alpha;\overline{\Omega}} + \|\rho\|_{1,\alpha;\overline{\Omega}}, \\
\text{where} \quad &C_a^{1,\alpha}(\overline{\Omega}) := \{w \in C^{1,\alpha}(\overline{\Omega}) : \int_{\Omega} w\, dx = 0\},
\end{aligned} \tag{2.12}$$

and $C^{m,\alpha}$ ($C^{\alpha} := C^{0,\alpha}$) denotes the usual Hölder spaces, for $0 < \alpha < 1$, with associated norms, cf. [8]. By virtue of (2.10)$_3$, $\mathbf{x} \mapsto \mathbf{A}(\lambda)\mathbf{x}$ is injective on $\partial\Omega$ for all $\lambda \in \mathbb{R}$, and any $\mathbf{f} \in C(\overline{\Omega}, \mathbb{R}^3) \cap C^1(\Omega, \mathbb{R}^3)$ satisfying (2.7)$_{2,3}$ ($\lambda$ fixed) is injective on $\overline{\Omega}$, cf. [3]. Hence, the integral of (2.7)$_2$ over the domain is zero, and likewise for (2.9)$_2$ (given the definition of $\mathbf{u}$). System (2.9) can then be expressed abstractly via

$$F(\lambda, w) = 0, \tag{2.13}$$

where $F : \mathbb{R} \times X \to Y$ is a $C^1$ mapping by virtue of (2.10), cf. [25].

**Remark 2.1.** The choice (2.7)$_3$ is made here for convenience. Our results that follow are easily generalized to the case $\mathbf{f}(\mathbf{x}) = \mathbf{g}(\lambda, \mathbf{x})$ $\mathbf{x} \in \partial\Omega$, for $\mathbf{g}(\cdot)$ sufficiently smooth, $\mathbf{g}(0, \mathbf{x}) \equiv \mathbf{x}$, and such that for each $\lambda \in \mathbb{R}$,

$$\mathbf{x} \mapsto \mathbf{g}(\lambda, \mathbf{x}) \text{ is injective on } \bar{\Omega},$$
$$\text{with } \int_\Omega \det \nabla_x \mathbf{g} \, dx = |\Omega|. \tag{2.14}$$

Condition (2.14)$_2$ insures global-volume compatibility. In this case, $\nabla_x \mathbf{g}(\lambda, x)$ replaces $\mathbf{A}(\lambda)$ in (2.9)$_{1,2}$, while an additional, additive term arises in (2.9)$_1$, which can be "absorbed" into the forcing term $\mathbf{b}$, while still satisfying (2.11)$_1$. Of course conditions (2.14) are trivially satisfied by the special case (2.7)$_3$, cf. (2.10)$_3$, (2.11)$_2$.

## 3. Local Solution Path

By virtue of (2.2), (2.4), (2.9) and (2.11), we verify that
$$F(0,0) = 0. \tag{3.1}$$

Our goal here is to employ the implicit function theorem to establish the existence of a local path solutions in a small neighborhood of the known solution point $(0,0) \in \mathbb{R} \times X$. Accordingly, we examine the Fréchet derivative of $w \mapsto F(\lambda, w)$ at $(0,0)$, denoted $T_o := D_w F(0,0) \in L(X,Y)$, given by

$$T_o[\mathbf{h}, r] = (A_o[\mathbf{h}, r], B_o[\mathbf{h}]), \tag{3.2}$$

where
$$A_o[\mathbf{h}, r] := \nabla \cdot (\mathcal{C}(\mathbf{I})[\nabla \mathbf{h}]) - \nabla r \text{ in } \Omega,$$
$$B_o[\mathbf{h}] := \nabla \cdot \mathbf{h} \text{ in } \Omega. \tag{3.3}$$

The following is well known:

**Proposition 3.1.** $T_o \in L(X,Y)$ *is injective.*

**Proof.** Consider $T_o[\mathbf{h}, r] = 0$. We take the dot product of (3.3)$_1$ with $\mathbf{h}$ and integrate over the domain. Integration by parts, while employing (3.3)$_2$ set equal to zero and $\mathbf{h} = \mathbf{0}$ on $\partial\Omega$, yields

$$\int_\Omega \nabla \mathbf{h} \cdot (\mathcal{C}(\mathbf{I})[\nabla \mathbf{h}]) dx = 0. \tag{3.4}$$

From here the argument is essentially the same as that given in [15; p. 37], viz., by taking the Fourier transform and exploiting strong ellipticity (2.6), which we do not repeat. One caveat here: Let $\hat{\mathbf{h}}(\mathbf{x})$ denote the Fourier transform of $\mathbf{h}$ (after extending $\mathbf{h}$ to zero on $\mathbb{R}^3 - \Omega$). The Fourier transform of $\nabla \mathbf{h}$ is then found to be $-i\hat{\mathbf{h}}(\mathbf{x}) \otimes \mathbf{x}$, and $\nabla \cdot \mathbf{h} = 0 \Rightarrow \hat{\mathbf{h}}(\mathbf{x}) \cdot \mathbf{x} = 0$ on $\mathbb{R}^3$, which is in consonance with (2.6) at $\mathbf{F} = \mathbf{I}$. The argument given in [15] then implies $\mathbf{h} = \mathbf{0}$ in $\bar{\Omega}$. Hence, $r \equiv const. = 0$, cf. (2.12). □

**Proposition 3.2.** $T_o \in L(X,Y)$ *is surjective.*

**Proof.** Let $\mathcal{I} : L(\mathbb{R}^3) \to L(\mathbb{R}^3)$ denote the identity, viz., in Cartesian components, $\mathcal{I}_{ijkl} := \delta_{ik}\delta_{jl}$, and define the one-parameter family of fourth-order tensors

$$\mathcal{C}_\mu := \mu\mathcal{I} + (1-\mu)\mathcal{C}(\mathbf{I}), \ \mu \in [0,1]. \tag{3.5}$$

Clearly $\mathcal{C}_\mu$ satisfies strong ellipticity (2.6) (with $\mathbf{F} = \mathbf{I}$) uniformly for all $\mu \in [0,1]$. Now consider the boundary value problem

$$A_\mu[\mathbf{h}, r] := \nabla \cdot (\mathcal{C}_\mu[\nabla \mathbf{h}]) - \nabla r = \boldsymbol{\tau} \text{ in } \Omega,$$
$$B_o[\mathbf{h}] := \nabla \cdot \mathbf{h} = \rho \text{ in } \Omega, \tag{3.6}$$
$$\mathbf{h} = \mathbf{0} \text{ on } \partial\Omega,$$

with $(\boldsymbol{\tau}, \rho) \in Y$, cf. (2.12). We also define

$$T_\mu := (A_\mu, B_o) \in L(X, Y). \tag{3.7}$$

Observe that (3.6) reduces to a Stokes problem at $\mu = 1$, viz., $(3.6)_1$ reads $A_1[\mathbf{h}, r] := \Delta \mathbf{h} - \nabla r = \boldsymbol{\tau}$ on $\Omega$; it is well-known that

$$T_1[\mathbf{h}, r] = (\boldsymbol{\tau}, \rho) \tag{3.8}$$

has a unique solution $(\mathbf{h}, r) \in X$ for all $(\boldsymbol{\tau}, \rho) \in Y$, cf. [24]. In particular, we conclude that $T_1$ has Fredholm index zero, where the Fredholm index of $T_\mu$ is defined by dim *Null* $T_\mu$ − dim *Range* $T_\mu$.

Next we claim, as a consequence of (2.6), that $\mathcal{C}_\mu$ satisfies the ellipticity condition of [1], which here reads

$$\det \begin{bmatrix} \mathbf{Q}_\mu(\mathbf{m}) & -\mathbf{m} \\ \mathbf{m} \cdot (\cdot) & 0 \end{bmatrix} \neq 0 \text{ for all } \mathbf{m} \in S^2, \tag{3.9}$$

where $\mathbf{Q}_\mu(\mathbf{m}) \in L(\mathbb{R}^3)$ is the acoustic tensor defined by $\mathbf{Q}_\mu(\mathbf{m})\mathbf{a} := \mathcal{C}_\mu[\mathbf{a} \otimes \mathbf{m}]\mathbf{m}$ for all $\mathbf{a} \in \mathbb{R}^3$, cf. [28]. The connection with the notation of [1] is readily established: In components, equations $(3.6)_{1,2}$ take the form

$$\sum_{j=1}^{4} \ell_{ij}(\partial) v_j = f_i, \ i = 1,...,4. \tag{3.10}$$

where $\mathbf{v} = (\mathbf{h}, r)$, $\mathbf{f} = (\boldsymbol{\tau}, \rho)$, and

$$\ell_{ij}(\mathbf{m}) = \sum_{k,l=1}^{3} \mathcal{C}_{\mu, ijkl} m_k m_l, \ i, j = 1, 2, 3,$$
$$\ell_{4j}(\mathbf{m}) = -\ell_{j4}(\mathbf{m}) = m_j, \ j = 1, 2, 3, \tag{3.11}$$
$$\ell_{44}(\mathbf{m}) = 0.$$

Choosing "weights" $s_i = 0$, $t_i = 2$, $i = 1, 2, 3$, $s_4 = -1$, $t_4 = 1$, we then verify $\deg \ell_{ij}(\mathbf{m}) \leq s_i + t_j$, $i, j = 1, ..., 4$. The components of the $4 \times 4$ matrix involved in (3.9) are given in (3.11). The former has a trivial null space by virtue of strong ellipticity (2.6), as we now demonstrate:

$$\begin{bmatrix} \mathbf{Q}_\mu(\mathbf{m}) & -\mathbf{m} \\ \mathbf{m} \cdot (\cdot) & 0 \end{bmatrix} \begin{pmatrix} \mathbf{a} \\ \zeta \end{pmatrix} = \begin{pmatrix} \mathbf{0} \\ 0 \end{pmatrix} \Rightarrow$$

$\mathbf{Q}_\mu(\mathbf{m})\mathbf{a} = \zeta\mathbf{m}$ and $\mathbf{m}\cdot\mathbf{a} = 0$. Hence, $0 = \mathbf{a}\cdot(\mathbf{Q}_\mu(\mathbf{m})\mathbf{a}) = (\mathbf{a}\otimes\mathbf{m})\cdot(\mathcal{C}_\mu[\mathbf{a}\otimes\mathbf{m}])$. In view of (2.6) (at $\mathbf{F} = \mathbf{I}$), we conclude that $\mathbf{a} = \mathbf{0}$, $\zeta = 0$ is the only solution, i.e., (3.9) is satisfied.

We note that the complementing condition is always satisfied in the presence of strong ellipticity (2.6) and Dirichlet boundary conditions (3.6)$_3$, cf. [1; p. 43-44]. Because (3.9) is also fulfilled, we then have the a-priori estimate

$$\|(\mathbf{h},r)\|_X \leq C\left[\|T_\mu[\mathbf{h},r]\|_Y + \|(\mathbf{h},r)\|_\infty\right] \quad \forall\ (\mathbf{h},r) \in X, \mu \in [0,1], \tag{3.12}$$

where $\|\cdot\|_\infty$ denotes the maximum norm over $\overline{\Omega}$, and where the constant $C$ is independent of $\mu$ and $(\mathbf{h},r)$, cf. [1]. This, in turn, implies that $T_\mu$ has a finite-dimensional kernel and a closed range, cf. [19], [26, p. 180]. The estimate (3.12) also implies that the Fredholm index of $T_\mu$ is constant on $[0,1]$ (cf. [13]), which gives the desired result, i.e., the injective map $T_o$ has Fredholm index zero. □

In view of Propositions 3.1 and 3.2, the implicit function theorem yields

**Proposition 3.3.** *There exists a unique local branch of solutions of* (2.13) *of the form*

$$\{(\lambda,\tilde{w}(\lambda)) \in \mathbb{R}\times X : |\lambda| < \varepsilon\},\ F(\lambda,\tilde{w}(\lambda)) \equiv 0, \tag{3.13}$$

*where $\varepsilon > 0$ is sufficiently small, $\lambda \mapsto \tilde{w}(\lambda)$ is $C^1$, and $\tilde{w}(0) = 0$. In particular, $w = 0$ is an isolated solution of $F(0,w) = 0$ in a sufficiently small neighborhood of $0 \in X$.*

**Remark 3.4.** If the incompressible, hyperelastic body is inhomogeneous, i.e., if the stored-energy density $W$ depends smoothly on $\mathbf{x} \in \overline{\Omega}$, then strong ellipticity (2.6) at $\mathbf{F} = \mathbf{I}$ does not imply uniqueness, i.e., Proposition 3.1 is no longer valid. However, if we make the physically reasonable assumption in this case that $\mathcal{C}(\mathbf{I},\mathbf{x})$ is positive definite on elements of $TS(\mathbb{R}^3)$, then uniqueness holds (cf. [15]), and the results of this section (and all that follow) remain true.

## 4. Global Continuation

We prove in this section that the local curve (3.13) is part of a global branch of solutions, i.e., a connected, locally compact set of solution pairs of (2.13). We now state this precisely:

**Theorem 4.1.** *Let $\mathbb{S} \subset \mathbb{R}\times X$ denote the solution set of* (2.13), *and let $\mathcal{C} \subset \mathbb{R}\times X$ be the connected component of $\mathbb{S}$ that contains the local solution curve* (3.13). *Then there are $\mathcal{C}^+,\mathcal{C}^- \subset \mathcal{C}$ such that*

(i) $\mathcal{C} = \mathcal{C}^+ \cup \{(0,0)\} \cup \mathcal{C}^-, \mathcal{C}^+ \cap \mathcal{C}^- = \varnothing$, *with $\mathcal{C}^+$ and $\mathcal{C}^-$ each unbounded in $\mathbb{R}\times X$,*

   *or*

(ii) $\mathcal{C}\setminus\{(0,0)\}$ *is connected.*

*In addition, for any $(\lambda,w) = (\lambda,\mathbf{u},p) \in \mathcal{C}$, the associated deformation $\mathbf{f}(\mathbf{x}) = \mathbf{A}(\lambda)\mathbf{x} + \mathbf{u}(\mathbf{x})$ is injective on $\overline{\Omega}$.*

**Remarks 4.2.** The solution sets $\mathcal{C}^+$ and $\mathcal{C}^-$ are defined in the proof of Theorem 4.1, cf. (4.29). The global solution branch $\mathcal{C}$ may or may not be unbounded in case (ii).

This result is more or less well known in cases for which the Leray-Schauder degree is applicable, i.e., when the nonlinear map can be converted to a compact perturbation of the identity. We refer to [23], [2] and in particular, to [14], where the precise degree-theoretic arguments yielding this version of the theorem are provided. The additional claim of injectivity follows as discussed after (2.12). We first verify the hypotheses of the nonlinear Fredholm degree of [5], [21], in order to deduce our theorem. We carry this out via two lemmas that follow.

Let $T(\lambda, w) := D_w F(\lambda, w) \in L(X, Y)$ denote the Fréchet derivative at $(\lambda, w) \in \mathbb{R} \times X$, given by

$$T(\lambda, \mathbf{u})[\mathbf{h}, r] = (L(\lambda, \mathbf{u})[\mathbf{h}, r], B(\lambda, \mathbf{u})[\mathbf{h}]) \quad \forall (\mathbf{h}, r) \in X, \tag{4.1}$$

where the components of principal part of the operator have the form

$$\begin{aligned} L(\lambda, \mathbf{u})[\mathbf{h}, r] &= \mathcal{C}(\mathbf{A}(\lambda) + \nabla \mathbf{u})\nabla^2 \mathbf{h} - Cof(\mathbf{A}(\lambda) + \nabla \mathbf{u})\nabla r + \ldots \text{ in } \Omega, \\ B(\lambda, \mathbf{u})[\mathbf{h}] &= Cof(\mathbf{A}(\lambda) + \nabla \mathbf{u}) \cdot \nabla \mathbf{h} \text{ in } \Omega; \end{aligned} \tag{4.2}$$

the latter follows from the chain rule and the fact that

$$\frac{d}{d\mathbf{F}} \det \mathbf{F} = Cof\, \mathbf{F}. \tag{4.3}$$

**Lemma 4.3.** *For each* $(\lambda, w) \in \mathbb{R} \times X$, $T(\lambda, w) \in L(X, Y)$ *is a Fredholm operator of index zero.*

**Proof.** The proof is similar to that of Proposition 3.2. For a given $(\lambda, w) \in \mathbb{R} \times X$, the ellipticity condition of [1] reads

$$\det \begin{bmatrix} \mathbf{Q}(\mathbf{F}_\lambda; \mathbf{m}) & -\hat{\mathbf{m}} \\ \hat{\mathbf{m}} \cdot (\cdot) & 0 \end{bmatrix} \neq 0 \text{ for all } \mathbf{m} \in S^2, \, \mathbf{F}_\lambda \in U(\mathbb{R}^3), \tag{4.4}$$

where

$$\begin{aligned} \mathbf{F}_\lambda &= \mathbf{A}(\lambda) + \nabla \mathbf{u}, \, (\mathbf{u}, p) \in \mathcal{B}_M, \\ \mathbf{Q}(\mathbf{F}_\lambda; \mathbf{m}) &\in L(\mathbb{R}^3), \, \hat{\mathbf{m}} := [Cof \mathbf{F}_\lambda]\mathbf{m}, \\ \mathbf{Q}(\mathbf{F}_\lambda; \mathbf{m})\mathbf{a} &= \mathcal{C}(\mathbf{F}_\lambda)[\mathbf{a} \otimes \mathbf{m}]\mathbf{m} \text{ for all } \mathbf{a} \in \mathbb{R}^3, \end{aligned} \tag{4.5}$$

cf. [28]. The connection with the notation of [1] is similar to that of (3.10), (3.11), viz., (4.2) corresponds to

$$\sum_{j=1}^{4} \ell_{ij}(\mathbf{m})(\partial) v_j = f_i, \, i=1,\ldots,4,$$

where $\mathbf{v} = (\mathbf{h}, r)$, $\mathbf{f} = (\boldsymbol{\tau}, \rho)$, and

$$\begin{aligned} \ell_{ij}(\mathbf{m}) &= \sum_{k,l=1}^{3} \mathcal{C}(\mathbf{F}_\lambda)_{ijkl} m_k m_l, \, i, j = 1, 2, 3, \\ \ell_{4j}(\mathbf{m}) &= -\ell_{j4}(\mathbf{m}) = \hat{m}_j, \, j = 1, 2, 3, \\ \ell_{44}(\mathbf{m}) &= 0, \end{aligned}$$

which give the components of the $4\times 4$ matrix in (4.4). We choose the same weights employed in (3.11).

By virtue of strong ellipticity (2.6), we claim that the $4\times 4$ matrix in (4.4) has a trivial null space:

$$\begin{bmatrix} \mathbf{Q}(\mathbf{F}_\lambda;\mathbf{m}) & -\hat{\mathbf{m}} \\ \hat{\mathbf{m}}\cdot(\cdot) & 0 \end{bmatrix}\begin{pmatrix} \mathbf{a} \\ \zeta \end{pmatrix} = \begin{pmatrix} \mathbf{0} \\ 0 \end{pmatrix} \Rightarrow$$

$\mathbf{Q}(\mathbf{F}_\lambda;\mathbf{m})\mathbf{a} = \zeta\hat{\mathbf{m}}$ and $[Cof\mathbf{F}_\lambda]\mathbf{m}\cdot\mathbf{a} = \hat{\mathbf{m}}\cdot\mathbf{a} = 0$. Hence, $0 = \mathbf{a}\cdot(\mathbf{Q}(\mathbf{F}_\lambda;\mathbf{m})\mathbf{a}) = (\mathbf{a}\otimes\mathbf{m})\cdot(\mathcal{C}(\mathbf{F}_\lambda)[\mathbf{a}\otimes\mathbf{m}])$. In view of (2.6), we conclude that $\mathbf{a} = \mathbf{0}$, $\zeta = 0$ is the unique solution, i.e., (4.4) is satisfied.

With this in hand, and again recalling that the complementing condition is satisfied (given strong ellipticity (2.6) with Dirichlet boundary conditions), we have the uniform a-priori estimate:

$$\|(\mathbf{h},r)\|_X \leq C\big[\|T(\lambda,\mathbf{u})[\mathbf{h},r]\|_Y + \|(\mathbf{h},r)\|_\infty\big] \quad \forall\ (\mathbf{h},r)\in X, \tag{4.6}$$

where the constant $C > 0$ is independent of $(\mathbf{h},r)$, cf. [1]. As before, we then deduce that $T(\lambda,w)\in L(X,Y)$ possesses a finite-dimensional kernel and a closed range. Moreover, with (4.6) in hand, it follows that the Fredholm index of $T(\lambda,w)$ is constant on the connected set $\mathbb{R}\times X$, cf. [13]. In particular, $T(0,0) \equiv T_o$ has Fredholm index zero, cf. Propositions 3.1, 3.2. □

**Lemma 4.4.** *The mapping $F$ is proper on bounded subsets, i.e., $F^{-1}(K)\cap \bar{D}$ is compact for each bounded set $D\subset \mathbb{R}\times X$ and compact set $K\subset Y$.*

**Proof.** We express the first component of (2.13), given explicitly by (2.9)$_1$, in convenient operator form via

$$F_1(\lambda,w) = \tilde{L}(\lambda,\mathbf{u})[\mathbf{u},p] + \varphi(\lambda,\mathbf{u}), \tag{4.7}$$

where

$$\begin{aligned}\tilde{L}(\lambda,\mathbf{u})[\mathbf{u},p] &:= \mathcal{C}(\mathbf{A}(\lambda)+\nabla\mathbf{u})\nabla^2\mathbf{u} - Cof(\mathbf{A}(\lambda)+\nabla\mathbf{u})\nabla p,\\ \varphi(\lambda,\mathbf{u}) &:= \mathbf{b}(\lambda,\nabla\mathbf{u},\mathbf{u},\cdot).\end{aligned} \tag{4.8}$$

Let $\{(\boldsymbol{\tau}_j,\rho_j)\}\subset K$ denote a convergent sequence in $Y$, and let the bounded sequence $\{(\lambda_j,\mathbf{u}_j,p_j)\}\subset \bar{D}\subset \mathbb{R}\times X$ satisfy

$$F(\lambda_j,\mathbf{u}_j,p_j) = (\boldsymbol{\tau}_j,\rho_j),\ j=1,2,... \tag{4.9}$$

Our goal is to show that $\{(\lambda_j,\mathbf{u}_j,p_j)\}$ has a convergent subsequence.

Since $D\subset \mathbb{R}\times X$ is bounded, we know from compact embedding, that possibly for a subsequence (not relabeled) we have, say,

$$\lambda_j \to \lambda\ \text{in}\ \mathbb{R},\ \mathbf{u}_j \to \mathbf{u}\ \text{in}\ C^2(\bar{\Omega},\mathbb{R}^3),\ \text{and}\ p_j \to p\ \text{in}\ C^1(\bar{\Omega},\mathbb{R}^3). \tag{4.10}$$

From (2.10) we then deduce that

$$\mathbf{A}(\lambda_j) \to \mathbf{A}(\lambda),$$
$$\varphi(\lambda_j, \mathbf{u}_j) \to \varphi(\lambda, \mathbf{u}) \text{ in } C^\alpha(\overline{\Omega}, \mathbb{R}^3). \tag{4.11}$$

In addition, the coefficients of the quasilinear operator in (4.8)$_1$, $\tilde{L}(\lambda_j, \mathbf{u}_j)$, are equicontinuous in $j \in \mathbb{N}$. Hence,

$$\left(\tilde{L}(\lambda_j, \mathbf{u}_j) - \tilde{L}(\lambda_k, \mathbf{u}_k)\right)[\mathbf{u}_j, p_j] \to 0 \text{ in } C^\alpha(\overline{\Omega}, \mathbb{R}^3) \text{ as } j, k \to \infty. \tag{4.12}$$

Next, for fixed $(\lambda, \mathbf{v}) \in \mathbb{R} \times C^{2,\alpha}(\overline{\Omega}, \mathbb{R}^3)$, we define a linear operator $\tilde{T}(\lambda, \mathbf{v}) \in L(X, Y)$ as follows:

$$\tilde{T}(\lambda, \mathbf{v})[\mathbf{h}, r] := (\tilde{L}(\lambda, \mathbf{v})[\mathbf{h}, r], B(\lambda, \mathbf{v})[\mathbf{h}]), \tag{4.13}$$

where $B(\lambda, \mathbf{v})$ is defined in (4.2)$_2$. On comparing (4.1), (4.2) with (4.8)$_1$, (4.13), we see that the linear operators $\tilde{T}(\lambda, \mathbf{u})$ and $T(\lambda, \mathbf{u})$ differ only by lower-order terms. Since $\{(\lambda_j, \mathbf{u}_j)\}$ is uniformly bounded, we then have, as in the proof of Lemma 4.2, the following uniform a-priori estimates:

$$\|(\mathbf{h}, r)\|_X \leq C\left[\|\tilde{T}(\lambda_j, \mathbf{u}_j)[\mathbf{h}, r]\|_Y + \|(\mathbf{h}, r)\|_\infty\right] \forall (\mathbf{h}, r) \in X,$$

where $C > 0$ is independent of $(\mathbf{h}, r)$ and the index $j \in \mathbb{N}$. In particular, we have

$$\|(\mathbf{u}_k, p_k) - (\mathbf{u}_\ell, p_\ell)\|_X \leq C\left[\|\tilde{T}(\lambda_j, \mathbf{u}_j)[\mathbf{u}_k - \mathbf{u}_\ell, p_k - p_\ell]\|_Y + \|(\mathbf{u}_k - \mathbf{u}_\ell, p_k - p_\ell)\|_\infty\right]. \tag{4.15}$$

In view of (4.10), the second term on the right side of inequality (4.15) approaches zero as $k, \ell \to \infty$.

We now make use of (2.9) and (4.9) to express the first term on the right side of (4.15) in the following equivalent form:

$$\begin{aligned}
\tilde{T}(\lambda_j, &\mathbf{u}_j)[\mathbf{u}_k - \mathbf{u}_\ell, p_k - p_\ell] \\
&= \tilde{T}(\lambda_j, \mathbf{u}_j)[\mathbf{u}_k - \mathbf{u}_\ell, p_k - p_\ell] \\
&\quad - F(\lambda_k, \mathbf{u}_k, p_k) + (\boldsymbol{\tau}_k, \rho_k) + F(\lambda_\ell, \mathbf{u}_\ell, p_\ell) - (\boldsymbol{\tau}_\ell, \rho_\ell) \\
&= \left(\tilde{L}(\lambda_j, \mathbf{u}_j)[\mathbf{u}_k - \mathbf{u}_\ell, p_k - p_\ell] - \tilde{L}(\lambda_k, \mathbf{u}_k)[\mathbf{u}_k, p_k] - \varphi(\lambda_k, \mathbf{u}_k) + \boldsymbol{\tau}_k \right. \\
&\quad + \tilde{L}(\lambda_\ell, \mathbf{u}_\ell)[\mathbf{u}_\ell, p_\ell] + \varphi(\lambda_\ell, \mathbf{u}_\ell) - \boldsymbol{\tau}_\ell, \; B(\lambda_j, \mathbf{u}_j)[\mathbf{u}_k - \mathbf{u}_\ell] \\
&\quad \left. - \det(\mathbf{A}(\lambda_k) + \nabla \mathbf{u}_k) + 1 + \rho_k + \det(\mathbf{A}(\lambda_\ell) + \nabla \mathbf{u}_\ell) - 1 - \rho_\ell\right).
\end{aligned} \tag{4.16}$$

Thus,

$$\begin{aligned}
\|\tilde{T}(\lambda_j, &\mathbf{u}_j)[\mathbf{u}_k - \mathbf{u}_\ell, p_k - p_\ell]\|_Y \\
&\leq \left\|\left(\tilde{L}(\lambda_j, \mathbf{u}_j) - \tilde{L}(\lambda_k, \mathbf{u}_k)\right)[\mathbf{u}_k, p_k]\right\|_{\alpha; \overline{\Omega}} + \left\|\left(\tilde{L}(\lambda_\ell, \mathbf{u}_\ell) - \tilde{L}(\lambda_j, \mathbf{u}_j)\right)[\mathbf{u}_\ell, p_\ell]\right\|_{\alpha; \overline{\Omega}} \\
&\quad + \|\varphi(\lambda_\ell, \mathbf{u}_\ell) - \varphi(\lambda_k, \mathbf{u}_k)\|_{\alpha; \overline{\Omega}} + \|\boldsymbol{\tau}_k - \boldsymbol{\tau}_\ell\|_{\alpha; \overline{\Omega}} + \|\rho_k - \rho_\ell\|_{1, \alpha; \overline{\Omega}} \\
&\quad + \|B(\mathbf{u}_j)[\mathbf{u}_k - \mathbf{u}_\ell] - \det(\mathbf{A}(\lambda_k) + \nabla \mathbf{u}_k) + \det(\mathbf{A}(\lambda_\ell) + \nabla \mathbf{u}_\ell)\|_{1, \alpha; \overline{\Omega}}.
\end{aligned} \tag{4.17}$$

We claim that the last term on the right side of (4.17) approaches zero as $j, k, \ell \to \infty$, which we prove below. Given this, then (4.10) - (4.12) and the convergence of $\{(\boldsymbol{\tau}_j, \rho_j)\}$ in $Y$ imply that the entire right side of (4.15) approaches zero as $j, k, \ell \to \infty$. Hence, inequality (4.15) shows that $\{(\mathbf{u}_k, p_k)\} \subset X$ is a Cauchy sequence, which completes the proof of the lemma.

To finish, we need to demonstrate that the last term on the right side of (4.17) goes to zero in the limit. First we note that (4.10) and continuity imply

$$\left\| B(\lambda_j, \mathbf{u}_j)[\mathbf{u}_k - \mathbf{u}_\ell] - \det(\mathbf{A}(\lambda_k) + \nabla \mathbf{u}_k) + \det(\mathbf{A}(\lambda_\ell) + \nabla \mathbf{u}_\ell) \right\|_{C^1(\bar{\Omega})} \to 0 \text{ as } j, k, \ell \to \infty. \quad (4.18)$$

So it's enough to show that the $C^\alpha$ norm of the gradient of the expression in (4.18) approaches zero as well. Define the fourth-order-tensor-valued function

$$\mathcal{D}(\mathbf{F}) := \frac{d}{d\mathbf{F}} Cof(\mathbf{F}), \; \mathcal{D}: GL^+(\mathbb{R}^3) \to L(L(\mathbb{R}^3)), \quad (4.19)$$

which is smooth. Although not needed here, we remark that an explicit but lengthy expression for $\mathcal{D}(\mathbf{F})$ can be obtained by first expressing $Cof\,\mathbf{F}$ in a convenient form via the Cayley-Hamilton theorem. In any case, (4.3), (4.19) and the chain rule lead to

$$\begin{aligned}
\nabla \Big( B(\lambda_j, \mathbf{u}_j)[\mathbf{u}_k - \mathbf{u}_\ell] &- \det(\mathbf{A}(\lambda_k) + \nabla \mathbf{u}_k) + \det(\mathbf{A}(\lambda_\ell) + \nabla \mathbf{u}_\ell) \Big) \cdot \mathbf{a} \\
&= \mathcal{D}(\mathbf{A}(\lambda_j) + \nabla \mathbf{u}_j)(\nabla^2 \mathbf{u}_j \mathbf{a}) \cdot (\nabla \mathbf{u}_k - \nabla \mathbf{u}_\ell) \\
&\quad + Cof(\mathbf{A}(\lambda_j) + \nabla \mathbf{u}_j) \cdot [(\nabla^2 \mathbf{u}_k - \nabla^2 \mathbf{u}_\ell) \mathbf{a}] \\
&\quad - Cof(\mathbf{A}(\lambda_k) + \nabla \mathbf{u}_k) \cdot (\nabla^2 \mathbf{u}_k \mathbf{a}) + Cof(\mathbf{A}(\lambda_\ell) + \nabla \mathbf{u}_\ell) \cdot (\nabla^2 \mathbf{u}_\ell \mathbf{a}) \\
&= \mathcal{D}(\mathbf{A}(\lambda_j) + \nabla \mathbf{u}_j)(\nabla^2 \mathbf{u}_j \mathbf{a}) \cdot (\nabla \mathbf{u}_k - \nabla \mathbf{u}_\ell) \\
&\quad + [Cof(\mathbf{A}(\lambda_j) + \nabla \mathbf{u}_j) - Cof(\mathbf{A}(\lambda_k) + \nabla \mathbf{u}_k)] \cdot (\nabla^2 \mathbf{u}_k \mathbf{a}) \\
&\quad + [Cof(\mathbf{A}(\lambda_\ell) + \nabla \mathbf{u}_\ell) - Cof(\mathbf{A}(\lambda_j) + \nabla \mathbf{u}_j)] \cdot (\nabla^2 \mathbf{u}_\ell \mathbf{a}) \text{ for all } \mathbf{a} \in \mathbb{R}^3.
\end{aligned} \quad (4.20)$$

Accordingly,

$$\begin{aligned}
\left\| \nabla \Big( B(\lambda_j, \mathbf{u}_j)[\mathbf{u}_k - \mathbf{u}_\ell] - \det(\mathbf{A}(\lambda_k) + \nabla \mathbf{u}_k) + \det(\mathbf{A}(\lambda_\ell) + \nabla \mathbf{u}_\ell) \Big) \right\|_{\alpha;\bar{\Omega}} \\
\leq \left\| \mathcal{D}(\mathbf{A}(\lambda_j) + \nabla \mathbf{u}_j) \nabla^2 \mathbf{u}_j (\nabla \mathbf{u}_k - \nabla \mathbf{u}_\ell) \right\|_{\alpha;\bar{\Omega}} \\
+ \left\| [Cof(\mathbf{A}(\lambda_j) + \nabla \mathbf{u}_j) - Cof(\mathbf{A}(\lambda_k) + \nabla \mathbf{u}_k)] \nabla^2 \mathbf{u}_k \right\|_{\alpha;\bar{\Omega}} \\
+ \left\| [Cof(\mathbf{A}(\lambda_\ell) + \nabla \mathbf{u}_\ell) - Cof(\mathbf{A}(\lambda_j) + \nabla \mathbf{u}_j)] \nabla^2 \mathbf{u}_\ell \right\|_{\alpha;\bar{\Omega}},
\end{aligned} \quad (4.21)$$

where the product algebras in (4.21) are understood as defined in (4.20). By the same reasoning used to obtain (4.12), viz., since $Cof(\mathbf{A}(\lambda_j) + \nabla \mathbf{u}_j)$ is equicontinuous in $j \in \mathbb{N}$, we find that the second and third terms on the right side of (4.21) each go to zero as $j, k, \ell \to \infty$. On the other hand, since $Cof(\cdot)$ is smooth on $U(\mathbb{R}^3)$, and $\{(\lambda_j, \mathbf{u}_j)\}$ is uniformly bounded, we have

$$\|\mathcal{D}(\mathbf{A}(\lambda_j)+\nabla\mathbf{u}_j)\nabla^2\mathbf{u}_j(\nabla\mathbf{u}_k-\nabla\mathbf{u}_\ell)\|_{\alpha;\bar{\Omega}}$$

$$\leq C\|\mathcal{D}(\mathbf{A}(\lambda_j)+\nabla\mathbf{u}_j)\|_{\alpha;\bar{\Omega}}\|\nabla^2\mathbf{u}_j\|_{\alpha;\bar{\Omega}}\|\nabla\mathbf{u}_k-\nabla\mathbf{u}_\ell\|_{\alpha;\bar{\Omega}} \quad (4.22)$$

$$\leq K\|\nabla^2\mathbf{u}_j\|_{\alpha;\bar{\Omega}}\|\nabla\mathbf{u}_k-\nabla\mathbf{u}_\ell\|_{\alpha;\bar{\Omega}},$$

for a constant $K>0$, independent of $j\in\mathbb{N}$. Finally (4.10) and the uniform boundedness of $\{\mathbf{u}_j\}$ in $C^{2,\alpha}(\bar{\Omega},\mathbb{R}^3)$ imply that the right side of (4.22) goes to zero as $j,k,\ell\to\infty$. □

We now sketch the construction of the Fredholm degree first introduced in [5] for $C^2$ maps and later generalized to the $C^1$ case in [21]. As mentioned in the Introduction, we refer to this as the FPR degree. In what follows, $GL(X,Y)\subset L(X,Y)$ denotes the set of all invertible linear operators, $\Phi_o(X,Y)\subset L(X,Y)$ is the set of all linear Fredholm operators of index zero, and $K(X)\subset L(X)$ denotes the set of all compact linear operators. The degree is based on the idea of the *parity* of a path of Fredholm operators with invertible endpoints. Specifically, with $I:=[0,1]$, let $T\in C^0(I,\Phi_o(X,Y))$ denote a path of Fredholm operators, and let $N\in C^0(I,GL(Y,X))$ be a parametrix, viz.,

$$N(t)T(t)=\mathcal{I}-\kappa(t),\ \forall t\in I, \quad (4.23)$$

where $\mathcal{I}$ denotes the identity on $X$ and $\kappa\in C^0(I,K(X))$. Assume that $T(0),T(1)\in GL(X,Y)$. Then the *parity* $\sigma(T,I)$ of $T$ on $I$ is defined by the product

$$\sigma(T,I):=ind(\mathcal{I}-\kappa(0))ind(\mathcal{I}-\kappa(1)), \quad (4.24)$$

where $ind(\mathcal{I}-\kappa)$ denotes the Leray-Schauder index. That is, $ind(\mathcal{I}-\kappa)=(-1)^d=\pm 1$, where $d=$ number of real eigenvalues of $\kappa\in K(X)$ strictly great than 1 (counted by algebraic multiplicity). We remark that parametrices always exists in this setting, while the parity (4.24) is independent of the choice of parametrix, cf. [5].

Next consider the nonlinear equation $G(w)=y\in Y$, with the restriction $G|_{\mathcal{O}}$, where $\mathcal{O}\subset X$ is simply connected, and let $\Upsilon\subset\mathcal{O}$ be open and bounded. Assume that $G\in C^1(\mathcal{O},Y)$, $G|_{\bar{\Upsilon}}$ is proper, and $DG(w)\in\Phi_0(X,Y)$ for all $w\in\mathcal{O}$, which is precisely our setting here (with $\mathcal{O}=X$ and $G(w)$ corresponding to $w\mapsto F(\lambda,w)$, cf. Lemmas 4.3 and 4.4). First, if there is no $p\in\mathcal{O}$ with $DG(p)\in GL(X,Y)$, then the absolute degree of $G$ on $\Upsilon$ is defined by

$$|d|(G,\Upsilon,y)=0. \quad (4.25)$$

If there is a $p\in\mathcal{O}$ with $DG(p)\in GL(X,Y)$, then $p$ is called a *base point*. Suppose that $y\notin G(\partial\Upsilon)$ is a regular value, i.e., either $\Upsilon\cap G^{-1}(y)=\emptyset$ or $DG(w)\in GL(X,Y)$ for every $w\in\Upsilon\cap G^{-1}(y)$. By properness and the inverse function theorem, it follows that, if not empty, then $\Upsilon\cap G^{-1}(y)=\{w_1,w_2,...,w_m\}$ for some finite number $m\in\mathbb{N}$, and the degree is defined by

$$d_p(G, \Upsilon, y) = \sum_{j=1}^{m} \sigma_j, \tag{4.26}$$

where $\sigma_j := \sigma(DG \circ \tau_j, I)$ is the parity, with $\tau_j$ being any continuous path on $I$ joining $p$ to $w_j$ in $\mathcal{O}$. In case $\Upsilon \cap G^{-1}(y) = \varnothing$, we define

$$d_p(G, \Upsilon, y) = 0. \tag{4.27}$$

It turns out that a change in base point produces the same degree (4.26) to within sign, the latter predicted by an explicit formula, cf. [5]. However, we have no need for that particular refinement here, since global continuation is readily deduced via the absolute degree, which is defined by (4.25) together with

$$|d|(G, \Upsilon, y) = |d_p(G, \Upsilon, y)|. \tag{4.28}$$

The singular-value case for the degree is handled in [5] via the Smale-Sard-Quinn theorem [22], which requires $G \in C^2$. An approximation theorem for $C^1$ Fredholm maps was later developed in [20] and subsequently employed in [21] to generalize the FPR degree to $G \in C^1$. In what follows, we employ the usual convention: $d_p(G, \Upsilon) := d_p(G, \Upsilon, 0)$ and $|d|(G, \Upsilon) := |d|(G, \Upsilon, 0)$.

The FPR degree possesses all of the usual properties of the Leray-Schauder degree, except that homotopy invariance no longer holds. Rather it is given to within sign – again, the latter predicted by a formula. Of course, the absolute degree (4.25), (4.28) is homotopy invariant but incapable of changing sign. While not appropriate for detecting bifurcation, this is all we require in order to complete the proof of Theorem 4.1.

**Proof of Theorem 4.1.** The proof is essentially that given in [14, p. 232], except that here we employ the FPR degree and the absolute degree. We first define

$$\begin{aligned} \mathcal{C}^+ &= \text{ component } \mathcal{P}_\varepsilon^+ \text{ in } \mathbb{S} \setminus \{(0,0)\}, \\ \mathcal{C}^- &= \text{ component } \mathcal{P}_\varepsilon^- \text{ in } \mathbb{S} \setminus \{(0,0)\}, \end{aligned} \tag{4.29}$$

where $\mathcal{P}_\varepsilon^{+(-)} := \{(\lambda, \tilde{w}(\lambda)) \in \mathbb{R} \times X : \lambda \in (0, \varepsilon)(\lambda \in (-\varepsilon, 0))\}$, for $\varepsilon > 0$ sufficiently small, are the "half paths" corresponding to Proposition 3.3. We note that the definitions of $\mathcal{C}^+$ and $\mathcal{C}^-$ are independent of $\varepsilon$ in (4.29). We also use the notation $\mathcal{P}_\varepsilon := \mathcal{P}_\varepsilon^+ \cup \{(0,0)\} \cup \mathcal{P}_\varepsilon^-$ for the local solution path. Now assume that $\mathcal{C}^+ \cap \mathcal{C}^- = \varnothing$, and suppose that $\mathcal{C}^+$ is bounded. We argue by contradiction to show that the latter is not possible.

We claim that there is a bounded open set $\Upsilon \subset \mathbb{R} \times X$ such that $\mathcal{C}^+ \subset \Upsilon$ and $\partial \Upsilon \cap \mathbb{S} = \{(0,0)\}$. First, by virtue of the implicit function theorem, there is a small open neighborhood of $(0,0)$, denoted $\mathcal{N}_o \subset \mathbb{R} \times X$, with $\mathcal{P}_\varepsilon \subset \mathcal{N}_o$ such that all solutions of (2.13) in $\mathcal{N}_o$ belong to $\mathcal{P}_\varepsilon$. Now define $\mathcal{N}_o^+ := \mathcal{N}_o \cap ((0, \infty) \times X)$. Taking $\varepsilon > 0$ sufficiently small, we have $\mathcal{P}_\varepsilon^+ \subset \mathcal{N}_o^+$, and we may assume that $(\lambda, w) \in \partial \mathcal{N}_o^+ \cap \mathbb{S}$ implies either $(\lambda, w) = (0,0)$ or $(\lambda, w) \in \mathcal{C}^+ \setminus \mathcal{P}_\varepsilon^+ := \mathcal{C}_\varepsilon^+$. Next, observe that $\mathcal{C}_\varepsilon^+$ is compact (by properness). By a well-known argument from [23], there is a bounded open set $\mathcal{U} \subset \mathbb{R} \times X$ such that $\mathcal{C}_\varepsilon^+ \subset \mathcal{U}$ and $\partial \mathcal{U} \cap (\mathbb{S} \setminus \mathcal{P}_\varepsilon^+) = \varnothing$. By construction, it follows that $\mathcal{N}_o^+ \not\subset U$; in particular, $(0,0) \notin U$. We provide the details pertaining to $U$ in the Appendix. In any case, $(\lambda, w) \in \partial \mathcal{U} \cap \mathbb{S} \Rightarrow (\lambda, w) \in \mathcal{P}_\varepsilon^+ \subset \mathcal{N}_o^+$. Finally we define $\Upsilon := \mathcal{U} \cup \mathcal{N}_o^+$, which has the desired properties: Clearly $\mathcal{C}^+ = \mathcal{C}_\varepsilon^+ \cup \mathcal{P}_\varepsilon^+ \subset \Upsilon$, and $\partial \Upsilon = \Gamma_1 \cup \Gamma_2$, where $\Gamma_1 := \partial \mathcal{N}_o^+ \setminus \mathcal{U}$ and $\Gamma_2 := \partial \mathcal{U} \setminus \mathcal{N}_o^+$. Then $\partial \Upsilon \cap \mathbb{S} = \Gamma_1 \cap \mathbb{S} = \{(0,0)\}$.

Let $B_\delta \subset X$ denote the open ball centered at the origin of radius $\delta > 0$. We choose $\delta$ sufficiently small such that $w = 0$ is the only solution of $F(0, w) = 0$ on $B_\delta$. We may assume that $\mathcal{N}_o \cap (\{0\} \times X) = \{0\} \times B_\delta$, and therefore $\{0\} \times B_\delta \subset \Gamma_1$. Define $\Upsilon_\lambda := \{w \in X : (\lambda, w) \in \Upsilon\}$ and $\overline{\Upsilon}_\lambda := \{w \in X : (\lambda, w) \in \overline{\Upsilon}\}$. For $\delta$ sufficiently small, we may assume that $\overline{(\Upsilon_0)} \cap \overline{B_\delta} = \varnothing$ and $\overline{(\Upsilon_0)} \cup \overline{B_\delta} = \overline{\Upsilon}_0$. In view of Propositions 3.1 and 3.2, $p = 0$ serves as a base point for $w \mapsto F(0, w)$, and by (4.28), the additivity of the base-point degree, and homotopy invariance of the absolute degree, we have

$$\left| d_0(F(0, \cdot), B_\delta) + d_0(F(0, \cdot), \Upsilon_0) \right| = |d|(F(\lambda, \cdot), \Upsilon_\lambda) \; \forall \lambda > 0, \qquad (4.30)$$
$$= |d|(F(\lambda, \cdot), \varnothing) = 0,$$

for $\lambda > 0$ sufficiently large. If $\Upsilon_0 = \varnothing$, then $d_0(F(\cdot, 0), \Upsilon_0) = 0$. On the other hand, if $\Upsilon_0 \neq \varnothing$, we also find

$$|d|(F(0, \cdot), \Upsilon_0) = |d|(F(\lambda, \cdot), \Upsilon_\lambda) \; \forall \lambda \leq 0, \qquad (4.31)$$
$$= |d|(F(\lambda, \cdot), \varnothing) = 0,$$

for $\lambda < 0$ sufficiently large in magnitude. Either way we conclude that $d_0(F(\cdot, 0), \Upsilon_0) = 0$, and then (4.30) implies $d_0(F(0, \cdot), B_\delta) = 0$. But this contradicts the fact that $w = 0$ is the only solution of $F(0, w) = 0$ on $B_\delta$: With base-point $p = 0$, let $\tau$ be any continuous path on $I$ connecting $w = 0$ to itself. In particular, $D_w F(0, \tau(0)) = D_w F(0, \tau(1)) = D_w F(0, 0)$, and (4.24) and (4.26) give

$$d_0(F(0, \cdot), B_\delta) = \sigma(D_w F(0, \cdot) \circ \tau, I) = 1.$$

We conclude that $\mathcal{C}^+$ is unbounded in $\mathbb{R} \times X$. Obviously we arrive at a similar conclusion when considering $\mathcal{C}^-$. In summary, we have two alternatives at this stage:

(i)' $\mathcal{C}^+ \cap \mathcal{C}^- = \varnothing$, with $\mathcal{C}^+$ and $\mathcal{C}^-$ each unbounded in $\mathbb{R} \times X$, or

(ii)' $\mathcal{C}^+ \cap \mathcal{C}^- \neq \varnothing$, i.e., $\mathcal{C}^+ = \mathcal{C}^-$.

It remains to show that $\mathcal{C} \setminus \{(0,0)\}$ has no other connected components distinct from $\mathcal{C}^+$ and $\mathcal{C}^-$, which is presumed in [14]. We argue by contradiction. Let $\mathcal{M}$ denote the union of all such components. Our goal is to demonstrate that $\mathcal{M} = \varnothing$. We observe first that $\mathcal{M} \cap \mathcal{N}_o = \varnothing$, where $\mathcal{N}_o$ is as defined above, courtesy of the implicit function theorem. We then make use of the following result from reference [4, p. 71, Prob. 6(b)]: Let $\mathcal{X}$ be a connected metric space containing at least two points. Let $\mathcal{Q}$ be a connected subset of $\mathcal{X}$ and $\mathcal{B}$ a connected component of $\mathcal{X} \setminus \mathcal{Q}$. Then $\mathcal{X} \setminus \mathcal{B}$ is connected. Case $(i)'$: First choose $\mathfrak{X} = \mathcal{C}$, $\mathcal{Q} = \{(0,0)\}$ and $\mathfrak{B} = \mathcal{C}^+$. Then $\mathcal{C} \setminus \mathcal{C}^+ = \mathcal{C}^- \cup \{(0,0)\} \cup \mathfrak{M}$ is connected. Next select $\mathfrak{X} = \mathcal{C}^- \cup \{(0,0)\} \cup \mathfrak{M}$ and $\mathfrak{B} = \mathcal{C}^-$, with $\mathcal{Q} = \{(0,0)\}$. Then $\{(0,0)\} \cup \mathfrak{M}$ is connected, which is a contradiction unless $\mathfrak{M} = \emptyset$. Alternative (*i*) of Theorem 4.1 then follows. Case (*ii*)′: Choose $\mathfrak{X} = \mathcal{C}$, $\mathcal{Q} = \{(0,0)\}$ and $\mathfrak{B} = \mathcal{C}^+$. Then $\mathcal{C} \setminus \mathcal{C}^+ = \{(0,0)\} \cup \mathfrak{M}$ is connected, and we again conclude that $\mathfrak{M} = \emptyset$. Thus, $\mathcal{C} \setminus \{(0,0)\} = \mathcal{C}^+ = \mathcal{C}^-$, which is the same as alternative (*ii*) of Theorem 4.1. □

## 5. Unbounded Solution Branches

We make two additional, physically reasonable hypotheses in this section that lead to the elimination of property (ii) of Theorem 4.1, and hence characterization (i) holds. First, we assume that the stored energy density has a global minimum at the identity, i.e.,

$$W(\mathbf{F}) > W(\mathbf{I}) = 0 \quad \forall \ \mathbf{F} \in U(\mathbb{R}^3) - SO(3). \tag{5.1}$$

With (5.1) in hand, it follows that

$$\int_D W(\mathbf{I} + \nabla \mathbf{v}) dx \geq 0, \tag{5.2}$$

for all bounded domains $D \subset \mathbb{R}^3$ and every $\mathbf{v} \in C^1(\overline{D}, \mathbb{R}^3)$ such that $\mathbf{v}|_{\partial D} = \mathbf{0}$ and $\det(\mathbf{I} + \nabla \mathbf{v}) = 1$ on $\overline{D}$, with equality in (5.2) holding only for $\mathbf{v} \equiv \mathbf{0}$. That is, $W(\cdot)$ is strictly quasiconvex at the identity. In addition, we assume that the domain $\Omega$ is star-shaped with respect to the origin, i.e.,

$$\mathbf{n}(\mathbf{x}) \cdot \mathbf{x} > 0 \quad \forall \ \mathbf{x} \in \partial \Omega. \tag{5.3}$$

**Theorem 5.1.** *Assume the previous hypotheses leading to Theorem 4.1 as well as* (5.1) *and* (5.3). *Then the global branch of solutions* $\mathcal{C}$ *of* (2.13) *is characterized solely by property* (i) *of Theorem 4.1. Moreover,* $\mathcal{C}^\pm \subset (0, \pm \infty) \times X$, *respectively.*

**Proof.** With (2.6), (5.2), (5.3) in hand, a result of [16] shows that system (2.7), (2.11) with homogeneous data ($\lambda = 0$) has the unique solution $\mathbf{u} = \mathbf{0}$, $p = const.$, the latter of which is zero by virtue of $(\mathbf{u}, p) \in X$, cf. (2.12)$_1$. Stated abstractly in terms of (2.13), we have $F(0, w) = 0 \Rightarrow w = 0$. In particular, there can be no solution $(0, w) \in \mathcal{C}$ with $w \neq 0$. Hence, $\mathcal{C} \setminus \{(0,0)\}$ is not connected. □

## 6. Concluding Remarks

Theorem 5.1 implies that there exists at least one solution for any given applied loading, and/or for some finite loading there are solutions of arbitrarily large norm, i.e., there exist solutions in the large. It is interesting to compare this with the same result obtained in [9] for the compressible case. There, based on the analogue of (5.1) and (5.3), the results of [16] also imply that the stress-free reference configuration corresponds to the unique solution of the problem in the absence of loading. Hence, the version of alternative (ii) in Theorem 4.1 cannot be true. However, this still leaves open the possibility that the "global" solution branch is bounded. This follows from the degree-theoretic construction, which is valid uniformly "away" from zero volume ratio. To eliminate the possibility that the branch is bounded, the stored energy function is hypothesized as the sum of two terms: One depends solely on the volume ratio, which grows unboundedly as its argument approaches zero, while the other term is presumed independent of the volume ratio. With this in hand, it is shown that the volume ratio can approach zero only in the unbounded limit along an unbounded solution branch. Here in the incompressible case, Theorem 5.1 follows directly without additional assumptions.

Hyperelasticity, assumed in Section 2, is not employed until Section 5. In fact the basic existence result in Theorem 4.1 requires only Cauchy elasticity with all other assumptions being the same. We emphasize that strict quasiconvexity, assumed in Section 5, is required only at the stress-free reference configuration. We also note that with the additional assumptions of Section 5 in hand, Theorem 5.1 can be proven directly without reliance on Theorem 4.1 (via the absolute degree [5], [21]) in a manner similar to that in [23].

In the case of "dead" loading, viz., the body force density $\mathbf{b}(\lambda, \cdot)$ is independent of the displacement field, a formulation with $W^{2,p}(\Omega, \mathbb{R}^3), W^{1,p}(\Omega)$ and $L^p(\Omega), p > 3,$ in place of $C^{2,\alpha}(\bar{\Omega}, \mathbb{R}^3), C^{1,\alpha}(\bar{\Omega})$ and $C^\alpha(\bar{\Omega})$, respectively (cf. (2.9), (2.12)) is readily carried out. Otherwise, with the same hypotheses, the results of Theorems 4.1 and 5.1 are valid in that setting. For "live" loadings, as considered here, it becomes necessary to work in higher-order Sobolev spaces to insure differentiability; we find Hölder spaces more convenient.

The absolute degree of [5], [21] is equally useful for global continuation in compressible nonlinear elasticity as treated in [10]. For example, this is carried out in [7] for compressible Mooney-Rivlin materials on unbounded domains. However, the FPR degree based on the parity can only be defined on simply-connected subsets $\mathcal{O} \subset X$ of a Banach space (here in this work $\mathcal{O} = X$). In realistic bifurcation problems of nonlinear elasticity with free-surface and/or traction conditions, the complementing condition plays a crucial role in the nonlinear Fredholm property. Given that the complementing condition can fail, the appropriate inequality conditions insuring its satisfaction are incorporated into the space of admissible deformations. One then naturally works in the connected component of that set containing the identity map. The situation is more complicated in the case of compressible nonlinear elasticity, given that the admissible set also contains only those deformations with strictly positive volume ratios. Either way, it is not clear that such a connected component of admissible solutions is simply-connected in an underlying Banach space. Indeed, the oriented degree of [10] was developed with a view toward global bifurcation problems in compressible nonlinear elasticity, cf. [12]. A systematic approach to global bifurcation in incompressible nonlinear elasticity awaits development.

While the proof of Theorem 4.1 closely follows the strategy employed in [14, p. 232] (via the Leray-Schauder degree), new features of general interest are provided here. First, a detailed construction of the open, bounded set $\Upsilon$ is given. In addition, no a-priori assumptions are made concerning the role of the solution continua $\mathcal{C}^\pm$ in the decomposition of the global solution branch $\mathcal{C}$. Rather, the precise decomposition is deduced after the global properties of $\mathcal{C}^\pm$ are established.

# Appendix

**Lemma A.1** *Let $\mathcal{C}^{\pm}$ and $\mathcal{P}_{\varepsilon}^{\pm}$ be as defined in* (4.29), *and set* $\mathcal{C}_{\varepsilon}^{+} = \mathcal{C}^{+} \setminus \mathcal{P}_{\varepsilon}^{+}$. *Suppose that* $\mathcal{C}^{+} \cap \mathcal{C}^{-} = \varnothing$, *with* $\mathcal{C}^{+} \subset \mathbb{R} \times X$ *bounded. Then there is a bounded open set* $\mathcal{U} \subset \mathbb{R} \times X$ *such that* $\mathcal{C}_{\varepsilon}^{+} \subset \mathcal{U}$, $\partial \mathcal{U} \cap (\mathcal{S} \setminus \mathcal{P}_{\varepsilon}^{+}) = \varnothing$, *and* $(0,0) \notin \mathcal{U}$.

**Proof.** Let $\mathcal{V}_{\delta}$ denote a $\delta$-neighborhood of the compact set $\mathcal{C}_{\varepsilon}^{+}$, for some $\delta$ satisfying $0 < \delta < \inf_{(\lambda,w) \in \mathcal{C}_{\varepsilon}^{+}} (|\lambda| + \|w\|_X)$. In particular, note that $(0,0) \notin \mathcal{V}_{\delta}$. Since $\mathcal{P}_{\varepsilon}^{+}$ is open in $\mathcal{S}$, it follows that $\mathcal{S} \setminus \mathcal{P}_{\varepsilon}^{+}$ is closed in $\mathcal{S}$. Define $\mathcal{K} = \overline{\mathcal{V}_{\delta}} \cap (\mathcal{S} \setminus \mathcal{P}_{\varepsilon}^{+})$. Then $\mathcal{K}$ is compact (by properness), and $\partial \mathcal{V}_{\delta} \cap \mathcal{C}_{\varepsilon}^{+} = \varnothing$. By virtue of the so-called Whyburn Lemma [27, Chpt. 1], there are disjoint compact subsets $\mathcal{K}_1, \mathcal{K}_2 \subset \mathcal{K}$ such that $\mathcal{C}_{\varepsilon}^{+} \subset \mathcal{K}_1$, $\partial \mathcal{V}_{\delta} \cap (\mathcal{S} \setminus \mathcal{P}_{\varepsilon}^{+}) \subset \mathcal{K}_2$, and $\mathcal{K} = \mathcal{K}_1 \cup \mathcal{K}_2$. Let $U$ be a sufficiently small open neighborhood of $\mathcal{K}_1$ such that $U \cap \mathcal{K}_2 = \varnothing$. Then $\mathcal{C}_{\varepsilon}^{+} \subset U$, and $\partial U \cap (\mathcal{S} \setminus \mathcal{P}_{\varepsilon}^{+}) = \varnothing$. □

*Acknowledgements*. This work was supported in part by the National Science Foundation through grant DMS-1613753, which is gratefully acknowledged. I thank Patrick Rabier for valuable suggestions on an earlier version of the manuscript. In particular, he showed me the last part of the proof of Theorem 4.1.

# References

1. Agmon, S., Douglis, A., Nirenberg, L.: Estimates near the boundary for solutions of elliptic partial differential equations satisfying general boundary conditions II, *Comm. Pure Appl. Math.*, **17,** 35–92, 1964.

2. Alexander, J.C., Yorke, J.A.: The implicit function theorem and global methods of cohomology, *J. Funct. Anal.*, **21,** 330–339, 1976.

3. Ciarlet, P.G.: *Mathematical Elasticity, Vol. I*, North-Holland, Amsterdam, 1988.

4. Dieudonné, J.: *Fouundations of Modern Analysis,* Academic Press, New York, 1969.

5. Fitzpatrick, P.M., Pejsachowicz, J., Rabier, P.J.: The degree of proper $C^2$ Fredholm mappings I, *J. Reine Angew. Math.*, **427**, 1–33, 1992.

6. Fosdick, R.L., MaSithigh, G.P.: Minimization in incompressible nonlinear elasticity theory, *J. Elasticity*, **16**, 267-301, 1986.

7. Gebran, H.G., Stuart, C.A.: Global continuation for quasilinear elliptic systems on $\mathbb{R}^n$ and the equations of elastostatics, *Adv. Nonlinear Stud.,* **9**, 727-762, 2009.

8. Gilbarg, D., Trudinger, N.S.: *Elliptic Partial Differential Equations of Second Order*, 2nd Ed., Springer-Verlag, New York, 1980.

9. Healey, T.J., Rosakis, P.: Unbounded branches of classical injective solutions in the forced displacement problem of nonlinear elastostatics, *J. Elasticity,* **49,** 65–78, 1997.

10. Healey, T.J., Simpson, H.C.: Global continuation in nonlinear elasticity, *Arch. Rational Mech. Anal.*, **143**, 1-28, 1998.

11. Healey, T.J.: Global continuation in displacement problems of nonlinear elastostatics via the Leray –

Schauder degree, *Arch. Rational Mech. Anal.*, **152**, 273-282, 2000.

12. Healey, T.J., Montes-Pizarro, E.L.: Global bifurcation in nonlinear elasticity with an application to barreling states of cylindrical columns, *J. Elasticity,* **71,** 33–58, 2003.

13. Kato, T.: *Perturbation Theory for Linear Operators*, 2nd Ed., Springer-Verlag, Berlin, 1980.

14. Kielhöfer, H.: *Bifurcation Theory*, 2nd Ed., Springer, New York, 2011.

15. Knops, R., Payne, L.: *Uniqueness Theorems in Linear Elasticity,* Springer-Verlag, New York, 1971.

16. Knops, R.J., Stuart, C.A.: Quasiconvexity and uniqueness of equilibrium solutions in nonlinear elasticity, *Arch. Rational Mech. Anal.*, **86**, 233-249, 1984.

17. Le Dret, H.: Constitutive laws and existence questions in incompressible nonlinear elasticity, *J. Elasticity* **15,** 369–387, 1985.

18. Leray, J., Schauder, J., Topologie et équations fonctionnelles, *Ann. Sci.Ec. Norm. Sup.* (3), **51**, 45–78, 1934.

19. Peetre, J.: Another approach to elliptic boundary problems*, Comm. Pure Appl. Math.,* **14,** 711–731, 1961.

20. Pejsachowicz, J., Rabier, P.J.: A substitute for the Sard-Smale theorem in the $C^1$ case, *J. Anal. Math.,* **76**, 289-319, 1998.

21. \_\_\_\_\_\_\_\_\_\_\_: Degree theory for $C^1$ Fredholm mappings of index 0, *J. Anal. Math.,* **76**, 265-288, 1998.

22. Quinn, F., Sard, A., Hausdorff conullity of critical images of Fredholm maps, *Amer. J. Math.,* **94,** 1101–1110, 1972.

23. Rabinowitz, P.H.: Some global results for nonlinear eigenvalue problems, *J. Funct. Anal.*, **7,** 487–513, 1971.

24. Solonnikov, V.A.: General boundary value problems for Douglis-Nirenberg Elliptic Systems II, in *Boundary value problems of mathematical physics, Part 4*, Trudy Mat. Inst. Steklov., **92**, 1966, 233–297; *Proc. Steklov Inst. Math.*, **92** (1968), 269–339.

25. Valent, T.: *Boundary Value Problems of Finite Elasticity*, Springer-Verlag, New York, 1988.

26. Wloka, J.: *Partial Differential Equations,* Cambridge University Press, Cambridge, 1987.

27. Whyburn, G.T.: *Topological Analysis,* Princeton University Press, Princeton, 1964.

28. Zee, L, Sternberg, E.: Ordinary and strong ellipticity in the equilibrium theory of incompressible hyperelastic solids, *Arch. Rational Mech. Anal.*, **83**, 53-90, 1983.